\documentclass[10pt,journal,a4paper,twoside]{IEEEtran}
\ifCLASSINFOpdf
\else
   \usepackage[dvips]{graphicx}
\fi
\usepackage{url}

\usepackage{graphicx}
\usepackage{amsmath}
\usepackage{amsfonts}
\usepackage{amssymb}
\usepackage[pdftex,dvipsnames]{xcolor}
\usepackage{multirow}
\usepackage{multicol}
\usepackage{enumerate}
\usepackage{tabularx}
\usepackage{bbm}
\usepackage{bm,amsthm}
\usepackage{hyperref}
\usepackage{comment}
\usepackage{enumitem}
\usepackage[labelformat=simple]{subcaption}
\usepackage{float}
\usepackage{blindtext}
\usepackage{balance}

\usepackage{dblfloatfix}% http://ctan.org/pkg/dblfloatfix
\usepackage{lipsum}% http://ctan.org/pkg/lipsum
\usepackage[backend=biber,style=ieee,doi=false,url=false,eprint=true, isbn = false]{biblatex}
\usepackage{xr}

\makeatletter
\newcommand*{\addFileDependency}[1]{% argument=file name and extension
  \typeout{(#1)}
  \@addtofilelist{#1}
  \IfFileExists{#1}{}{\typeout{No file #1.}}
}
\makeatother

\newcommand*{\myexternaldocument}[1]{%
    \externaldocument{#1}%
    \addFileDependency{#1.tex}%
    \addFileDependency{#1.aux}%
}

\myexternaldocument{IEEE_SP_Letter_Supplementary_Material_v2}

\newtheorem{lemma}{Lemma}

\newtheorem{remark}{Remark}
\newtheorem{prop}{Proposition}

\newtheorem{assumption}{AS}

\DeclareMathOperator*{\R}{\mathbb{R}}
\DeclareMathOperator*{\Z}{{\mathbb{Z}}^+}

\newcommand{\cf}{\emph{cf.}}
\newcommand{\tran}{^{\mbox{\scriptsize T}}}  % transpose

\usepackage{algorithm}% http://ctan.org/pkg/algorithm
\usepackage{algpseudocode}% http://ctan.org/pkg/algorithmicx

\usepackage{mathtools}

\newcommand{\addb}[1]{{{\color{black}#1}}}
\makeatletter
\newenvironment{breakablealgorithm}
  {% \begin{breakablealgorithm}
   \begin{center}
     \refstepcounter{algorithm}% New algorithm
     \hrule height.8pt depth0pt \kern2pt% \@fs@pre for \@fs@ruled
     \renewcommand{\caption}[2][\relax]{% Make a new \caption
       {\raggedright\textbf{\fname@algorithm~\thealgorithm} ##2\par}%
       \ifx\relax##1\relax % #1 is \relax
         \addcontentsline{loa}{algorithm}{\protect\numberline{\thealgorithm}##2}%
       \else % #1 is not \relax
         \addcontentsline{loa}{algorithm}{\protect\numberline{\thealgorithm}##1}%
       \fi
       \kern2pt\hrule\kern2pt
     }
  }{% \end{breakablealgorithm}
     \kern2pt\hrule\relax% \@fs@post for \@fs@ruled
   \end{center}
  }
\makeatother

\renewbibmacro{in:}{}
\DeclareBibliographyCategory{important}
\DeclareFieldFormat{labelprefix}{%
  \ifcategory{important}
    {\textcolor{blue!80!black}{#1}}
    {#1}}

\DeclareFieldFormat{labelnumber}{%
  \ifcategory{important}
    {\textcolor{blue!80!black}{#1}}
    {#1}}

\DeclareFieldFormat{extranumber}{%
  \ifcategory{important}
    {\textcolor{blue!80!black}{#1}}
    {#1}}

\DeclareFieldFormat{labelnumberwidth}{%
  \ifcategory{important}
    {\textcolor{blue!80!black}{\mkbibbrackets{#1}}}
    {\mkbibbrackets{#1}}}

\DeclareFieldFormat{article}{%
  \ifcategory{important}
    {\textcolor{blue!80!black}{\mkbibbrackets{#1}}}
    {\mkbibbrackets{#1}}}

\AtEveryBibitem{
\ifcategory{important}%
    {\color{blue!80!black}}%
    {}%
}

\begin{document}

\title{On the Convergence of Inexact Gradient Descent with Controlled Synchronization Steps}

\author{Sandushan Ranaweera, Chathuranga Weeraddana, Prathapasinghe Dharmawansa, and Carlo Fischione
%\thanks{}
%\thanks{}
\thanks{S. Ranaweera (e-mail: sandushan@ieee.org) and P. Dharmawansa (e-mail: prathapa@uom.lk) are with the Dept. of Electronic and Telecom. Eng., University of Moratuwa, Sri Lanka. 
C. Weeraddana (e-mail: Chathuranga.Weeraddana@oulu.fi) is with the Faculty of Information Technology and Electrical Eng., University of Oulu, Finland. C. Fischione (e-mail: carlofi@kth.se) is with the School of Electrical Eng. and Computer Science, KTH Royal Institute of Technology, Sweden.}}

% \markboth{IEEE Signal Processing Letters}
% {On the Convergence of Inexact Gradient Descent with Controlled Synchronization Steps}
%{Shell \MakeLowercase{\textit{et al.}}: Bare Demo of IEEEtran.cls for IEEE Journals}
\maketitle

%Towards Decentralized Distributed Gradient Descent with Imperfect Communications

\begin{abstract}

We develop a gradient-like algorithm to minimize a sum of peer objective functions based on coordination through a peer interconnection network. The coordination admits two stages: the first is to constitute a gradient, possibly with errors, for updating locally replicated decision variables at each peer and the second is used for error-free averaging for synchronizing local replicas. Unlike many related algorithms, the errors permitted in our algorithm can cover a wide range of inexactnesses, as long as they are bounded. \addb{Moreover, we do not impose any gradient boundedness conditions for the objective functions.}
Furthermore, the second stage is not conducted in a periodic manner, like many related algorithms. Instead, a locally verifiable criterion is devised to dynamically trigger the peer-to-peer coordination at the second stage, so that expensive communication overhead for error-free averaging can significantly be reduced. Finally, the convergence of the algorithm is established under mild conditions.

% The algorithm 
% consists of two stages performed in an iterative manner. The first is used to exchange gradients possibly with \emph{errors}. We have \emph{no} restrictions on the errors of local gradient estimates, except that they are \emph{bounded}. As a result, our modeling can handle errors beyond those of classic quantization models with restrictions, such as diminishing and unbiasedness. For instance, a cheap low-bit quantization can be used \emph{throughout} the algorithm iterates under the first stage. The second stage is used to \emph{error-free} avaraging of local replicas. 

% This article develops an algorithm to solve a decentralized optimization problem via gradient descent (GD) algorithm, where $N$ sub-systems (SS) collectively minimize a sum of strongly convex and L-smooth functions. In particular, we assign each SS a local function and they perform GD iterates by communicating local gradient information between SSs. Moreover, we allow additive and bounded errors in these gradient communications while compensating them with additional synchronization steps. Furthermore, we develop an adaptive algorithm to perform GD while choosing iterations to apply these synchronization steps based on locally available information at SSs. We show that with constant step size, the  proposed algorithm converges to a neighborhood of the optimal solution with a linear rate of convergence. Our simulation results show that 

%{\color{red} Who cares this??}{\color{blue}Finally, theoretical derivations are verified by numerical simulations.{\color{cyan} this part can removed}}
\end{abstract}

\begin{IEEEkeywords}
Distributed optimization, inexact algorithms

%{{\color{red} Do you really need this?}{\color{blue}}decentralized system}, {\color{red} I would put ``algorithm" and ``convex function" here} {\color{blue} only put first three index terms since the paper is focused om distributed optimization problem}
\end{IEEEkeywords}

\IEEEpeerreviewmaketitle

\section{Introduction}
\IEEEPARstart{G}radient descent and its variants often lend themselves fully amenable to parallel and distributed algorithms, which are highly desirable in {large-scale} optimization problems~\cite{bertsekas2015parallel_f}. As a result, solution methods for many problems of recent interest are predominantly based on such gradient-like algorithms~\cite{NIPS2012_6aca9700_f,Chilimbi2014ProjectSystem_f,BrendanMcMahan2017,Konecny2016FederatedIntelligence_f,Khaled2019,Stich2019LocalLittle_f,Wang2018_f,palomar_eldar_2009_f,Angelia2020DistributedOptimization_f,Xiao2004SimultaneousDecomposition_f,Low-Layering-as-Optimization-Decomposition-2007,Nedic2009DistributedOptimization_f,9479747_f}. \addb{Broadly speaking, those algorithms developments are twofold~\cite{Verbraeken2020_f}: a) a federated setting where a {central controller} (CC) intervenes for decision variable update~\cite{NIPS2012_6aca9700_f,Chilimbi2014ProjectSystem_f,BrendanMcMahan2017,Konecny2016FederatedIntelligence_f,Khaled2019,Stich2019LocalLittle_f, Wang2018_f}; b) a {peer-to-peer} (PP) setting where subsystems (SSs), each with its {replicated decision variable}, perform locally the update through some peer {interconnection network, often modeled by a connected graph}~\cite{palomar_eldar_2009_f,Angelia2020DistributedOptimization_f,Nedic2009DistributedOptimization_f,Xiao2004SimultaneousDecomposition_f,Low-Layering-as-Optimization-Decomposition-2007,9479747_f}. In this setting, the algorithm relies on neighbors specified by the graph and does not rely on a CC  like in the federated setting.} As such, it appears that PP setting is more appealing than the CC setting due to many reasons, such as higher scalability and inherently decentralized collection of big data sets, among others~\cite{bertsekas2015parallel_f,Verbraeken2020_f}.
%~\cite[\S1.1.1]{bertsekas2015parallel_f}~\cite[\S~3.4]{Verbraeken2020_f}
%
In the context of a PP setting, a more fundamental concern is that the distributed algorithms usually undergo {inevitable} inexact conditions, e.g., unreliable and often limited communication capabilities~\cite{bertsekas2015parallel_f,Verbraeken2020_f,Kairouz2021_f}. 
% In the context of a peer-to-peer setting, the predominant issues are clearly different from those of a CC setting with a central locus of control. A more fundamental concern is that the distributed algorithms under peer-to-peer setting usually undergo \emph{inevitable} inexact conditions, e.g., unreliable and often limited communication capabilities~\cite{bertsekas2015parallel_f,Verbraeken2020_f,Kairouz2021_f}.
Thus, unlike the inexactnesses under CC settings~\cite{Chen2022DistributedDifferences_f,sindri_2018,sindri_2020,Nedic2018NetworkOptimization_f,Magnusson2018CommunicationOptimization_f,9311407_f}, those under PP settings influence the optimality, convergence, and effective implementation of algorithms. Consequently, there is an appeal to design effective algorithms under PP setting~\cite{Lee2019FiniteConvergence_f,Nedic2008DistributedEffects_f,Koloskova2019DecentralizedCommunication_f,George2020DistributedCommunication_f,9296963_f,7462298_f,Alistarh2018TheMethods_f,Alistarh2017QSGD:Encoding_f,Yu2019ParallelLearning_f,Zhou2018OnOptimization_f,Xie2020CSER:Reset_f}.

%\cite{TuranRobustGradients_f,Xu2021DistributedDerivatives_f,Agarwal2012DistributedOptimization_f,Chen2018LAG:Learning_f,Li2021Communication-CensoredDescent_f,Chen2022DistributedDifferences_f} Kairouz2021_f

%\cite{}\cite{Kairouz2021_f}

Algorithms in \cite{Lee2019FiniteConvergence_f,Nedic2008DistributedEffects_f,Koloskova2019DecentralizedCommunication_f,George2020DistributedCommunication_f,9296963_f,7462298_f} are based on distributed subgradient methods due to \cite{Nedic2009DistributedOptimization_f}. Some of these methods consider quantization models~ \cite{Nedic2008DistributedEffects_f,Koloskova2019DecentralizedCommunication_f,Lee2019FiniteConvergence_f} and others consider event-triggered models~\cite{George2020DistributedCommunication_f,9296963_f,7462298_f}, so as to {reduce} the communication burden between SSs. 
% The convergence rates are usually on the order of $\mathcal{O}(1/k)$ or inferior, where $k$ is the iteration index of the algorithm. {\bf Authors in~\cite{Lee2019FiniteConvergence_f,7462298_f} have, however,} established linear convergence of related algorithms under more restrictive conditions.
The gradient boundedness of underlying objective functions, although a restriction, has been considered in~\cite{Lee2019FiniteConvergence_f,Nedic2008DistributedEffects_f,Koloskova2019DecentralizedCommunication_f,George2020DistributedCommunication_f,9296963_f,7462298_f}, a technical assumption that enables convergences. The errors introduced in~\cite{Lee2019FiniteConvergence_f,Nedic2008DistributedEffects_f,Koloskova2019DecentralizedCommunication_f,George2020DistributedCommunication_f,9296963_f,7462298_f} can be viewed as {controllable}, in the sense that they are at the disposal of the algorithm. For example, quantization models in~\cite{Koloskova2019DecentralizedCommunication_f,Lee2019FiniteConvergence_f} are chosen to be {unbiased}, a favorable condition for convergence. However, a peer interconnection network can often admit errors that are {not} at the disposal of the algorithm, e.g., wireless links~\cite[\S~9]{goldsmith_2005_f}, limiting the applicability of developments in~\cite{Lee2019FiniteConvergence_f,Nedic2008DistributedEffects_f,Koloskova2019DecentralizedCommunication_f,George2020DistributedCommunication_f,9296963_f,7462298_f}

\addb{Works in~\cite{Alistarh2018TheMethods_f,Alistarh2017QSGD:Encoding_f,Yu2019ParallelLearning_f,Zhou2018OnOptimization_f,Xie2020CSER:Reset_f} rely on PP coordination to constitute a gradient, in contrast to common federated settings where primal variables are coordinated instead. Then the resulting gradients are for updating their locally replicated decision variables.} They are persuaded again under quantization settings (e.g.,~\cite{Alistarh2018TheMethods_f,Alistarh2017QSGD:Encoding_f}) and event-triggered settings (e.g.,~\cite{Yu2019ParallelLearning_f,Zhou2018OnOptimization_f}). Hybrid variants have also been considered by some authors, e.g.,~\cite{Xie2020CSER:Reset_f}. 
%Convergence rates of the related algorithms are shown to be on the order of $\mathcal{O}(1/k)$. 
Similar to the developments noted in the preceding discussion, errors introduced in~\cite{Alistarh2018TheMethods_f,Alistarh2017QSGD:Encoding_f,Yu2019ParallelLearning_f,Zhou2018OnOptimization_f,Xie2020CSER:Reset_f} are also controlled by the algorithms. For example, the quantization models in~\cite{Alistarh2018TheMethods_f} and \cite{Alistarh2017QSGD:Encoding_f} are chosen so that the errors are diminishing and unbiased, respectively. Moreover, the authors in \cite{Yu2019ParallelLearning_f,Xie2020CSER:Reset_f} have specific impositions on the gradient boundedness. 

\addb{It is worth noting that many algorithms in either of the setting federated or PP (e.g., \cite{Khaled2019,Stich2019LocalLittle_f,Wang2018_f,Yu2019ParallelLearning_f,Zhou2018OnOptimization_f, ,Xie2020CSER:Reset_f}) have considered an averaging step performed at {periodic or predefined} epochs to enable the consistency of the locally replicates decision variables. Depending on the context, this entails periodic communication through the CC or through the PP interconnection network. From a communication overhead point of view, however, such an overhead for periodic communication seems like a restriction. This may be avoided by dynamically choosing the averaging epochs for synchronization.} 

% \addr{In particular, the subsystem coordination in our method is peer-to-peer that can be specified using a graph. Moreover, unlike common federated learning methods, gradients are communicated instead of primal decision variables. Thus, the proposed setting is directly applicable to dual decomposition methods \addr{[REF]}. \cite{Khaled2019,Stich2019LocalLittle_f}, event-triggered settings, errors introduced are also controlled by the underlying algorithm, ave specific impositions
% on the gradient boundednes, ave considered an averaging
% step performed at periodic or predefined epochs to enable the
% consistency of the locally replicates decision variables}

In this paper, we develop an algorithm that relies on PP coordination to constitute a gradient for updating locally replicated decision variables associated with a problem of minimizing the sum of peer objective functions. The algorithm iterates
two stages. The first is used to exchange gradients possibly with {errors}. We have {no} restrictions on the errors of local gradient estimates, except that they are {bounded}. As a result, our modeling can handle errors beyond those of classic quantization models with restrictions, such as diminishing and unbiasedness. For instance, a cheap low-bit quantization can be used {throughout} the algorithm iterates under the first stage. The second stage is used to {error-free} averaging for synchronizing local replicas. In this respect, unlike other related algorithms, we do {not} rely on periodic communication over the PP network. Instead, a {locally verifiable} criterion is devised to {dynamically} trigger the averaging step, only when necessary. This has the advantage of minimizing expensive communication overhead. \addb{Throughout this paper, we consider the PP network to be fully connected.~\footnote{\addb{An extension to an arbitrary graph is possible with an additional assumption on the gradient boundedness. The details are provided in the Appendix.}} Subsequently, the convergence of the algorithm is established and is shown to be {linear}.}

\section{Problem Formulation}
\label{sec:prob_formulation}

%Broadly speaking, we pose the problem: What if, in a multi subsystem (SS) setting, each SS performs \emph{locally} the classic gradient descent update of an \emph{own copy} of the decision variable with \emph{inexact gradient measurements}?

% {\color{red} (This is already outlined in the Introduction, I guess. Therefore, we should reduce this.) 
% Consider the scenario in which an each subsystem of a certain multi-subsystem setting evaluates locally the gradient descent update of {\color{red} its} \emph{own copy} of the decision variable with \emph{inexact gradient measurements}.
% %
% Since the local iterates are performed in parallel, there is no guarantee that at distinct SSs updates evolve in a meaningful manner with desirable convergence properties unless some control is imposed. Thus, the main problem in this research is to explore the convergence properties of such parallel iterates, possibly with specific modifications which are favourable from a practical standpoint. }
%
%{\color{red}How about reformulate like this? \bf Consider the scenario in which an each subsystem (SS) of a certain multi-subsystem evaluates its own gradient updates locally............. }

%In this section, we will formally define the optimization problem considered in this research.
%\subsection{Main Problem}
Consider $N$ peers or subsystems which solve the problem
\begin{equation} \label{eq:main-problem}
\begin{array}{ll}
\mbox{minimize} &  f(\boldsymbol{x}) = \sum_{i=1}^{N} f_{i} \left( \boldsymbol{x} \right)
\end{array}
\end{equation}
where $\boldsymbol{x} \in \mathbb{R}^{n}$ and $f_i:\R^n\rightarrow \R$, $ i\in \mathcal{N}\triangleq\{1,\ldots,N\}$, be a function satisfying the following standard assumption:
\begin{assumption} 
\label{Ass:objectives_assump}
The objective function $f_{i}$, $i\in\mathcal{N}$, is strongly convex with constant $\ell_{i}>0$ and is $L_i$-smooth, i.e., $\nabla f_{i}$ is Lipschitz continuous with the constant $L_{i}>0$.
\end{assumption} 
A commonly used iterative algorithms for solving problem~\eqref{eq:main-problem} is the gradient descent (GD) algorithm $ \boldsymbol{x}^{(k+1)}
    =\textstyle\boldsymbol{x}^{(k)} - \gamma \sum_{j=1}^{N} \nabla f_{j}\big(\boldsymbol{x}^{(k)}\big)$,
% \begin{align}
% \label{eq:gd_var_update-distributed}
%     \boldsymbol{x}^{(k+1)}
%     &=\textstyle\boldsymbol{x}^{(k)} - \gamma \sum_{j=1}^{N} \nabla f_{j}\big(\boldsymbol{x}^{(k)}\big) 
% \end{align}
where $k\in\Z\triangleq\{0,1,2,\ldots\}$ is the iteration index and $\gamma$ is the step size. In contrast, here we assume a setting where each subsystem (SS)~$i$ performs locally the variable GD update of its own copy ${\boldsymbol x}_i^{(k+1)}$ of ${\boldsymbol x}^{(k+1)}$. This setting facilitates a distributed implementation of GD and thus each SS~$i$ relies on a {communication} with SS~$j$ to get a rough measurement of $\nabla f_{j}\big(\boldsymbol{x}^{(k)}_j\big)$ as specified below:
% \begin{itemize}
%     \item \emph{Communication}: Each SS~$i$ communicates with other SSs $j\neq i$ to compute $\nabla f\big(\boldsymbol{x}^{(k)}_i\big)=\sum_{j=1}^{N} \nabla f_{j}\big(\boldsymbol{x}^{(k)}_j\big)$.
%      \item \emph{Computation}: Each SS~$i\in\mathcal{N}$ performs locally the variable update \eqref{eq:gd_var_update-distributed} of local copy ${\boldsymbol x}_i^{(k+1)}$ of ${\boldsymbol x}^{(k+1)}$ by using $\nabla f\big(\boldsymbol{x}^{(k)}_i\big)$.
% \end{itemize}
% 
% 
% 
% 
% 
% \subsection{Inexact Distributed Implementation}\label{subsec:Computation-and-Communication-Setting}
% Distributed implementation of \eqref{eq:gd_var_update-distributed} (i.e., GD) includes the following communication and computation stages.
% \begin{itemize}
%     \item \emph{Communication}: Each SS~$i$ communicates with other SSs $j\neq i$ to compute $\nabla f\big(\boldsymbol{x}^{(k)}_i\big)=\sum_{j=1}^{N} \nabla f_{j}\big(\boldsymbol{x}^{(k)}_j\big)$.
%      \item \emph{Computation}: Each SS~$i\in\mathcal{N}$ performs locally the variable update \eqref{eq:gd_var_update-distributed} of local copy ${\boldsymbol x}_i^{(k+1)}$ of ${\boldsymbol x}^{(k+1)}$ by using $\nabla f\big(\boldsymbol{x}^{(k)}_i\big)$.
% \end{itemize}
%{\color{red} In this respect, }The following non-ideal communication setting is assumed.
\begin{assumption} 
\label{Ass:Non-Ideal-Communication}
$\forall$ $i,j\in\mathcal{N}$, s.t. $i\neq j$, gradient measurement $\boldsymbol{h}_{ij}^{(k)}\in\R^n$ received by $i$-th SS from $j$-th SS at $k$-th iteration is given by 
\begin{equation}
\label{grad_est_model}
    \boldsymbol{h}_{ij}^{(k)} = \nabla f_{j} \big(\boldsymbol{x}_j^{(k)} \big) + \boldsymbol{\epsilon}_{ij}^{(k)}
\end{equation}
where $\boldsymbol{\epsilon}_{ij}^{(k)}\in\R^n$ is a error such that {$||\boldsymbol{\epsilon}_{ij}^{(k)}||\leq \epsilon$ with $||\cdot||$ denoting the Euclidean norm.}
\end{assumption} 
\addb{The parameters $\boldsymbol{\epsilon}_{ij}^{(k)}$ model measurement errors, noises, quantization errors\footnote{\cf~\cite[Definition~2]{sindri_2020} for such a quantization that yield an error as in \textbf{AS}~\ref{Ass:Non-Ideal-Communication}.} due to compression, among others. However, note the upper bound condition on $\boldsymbol{\epsilon}_{ij}^{(k)}$ in \textbf{AS}~\ref{Ass:Non-Ideal-Communication}, where $\epsilon$ can be thought of as the worst-case characteristic of
errors throughout the algorithm. Under \textbf{AS}~\ref{Ass:Non-Ideal-Communication}, the gradient $\nabla f\big(\boldsymbol{x}^{(k)}_i\big)$ is distorted, which in turn admits the following iterate:} % SS~$i$, $i\in\mathcal{N}$:}
\begin{equation}
\label{eq:local-iterates-inexact}
 \textstyle  \boldsymbol{x}_{i}^{(k+1)} = \boldsymbol{x}_{i}^{(k)} - \gamma  \sum_{j=1}^{N} \boldsymbol{h}^{(k)}_{ij}, \quad i\in\mathcal{N}.
\end{equation}
% \exists~k\in\Z\setminus~\{0\},~......given $\boldsymbol{x}^{(0)}_j=\boldsymbol{x}^{(0)}_i$ $\forall$ $i,j\in\mathcal{N}$.
Strictly speaking, the local variables updates should be {consistent} in the sense that $\forall~k\in\Z$, $\forall$ $i,j\in\mathcal{N}$, $\boldsymbol{x}^{(k)}_j=\boldsymbol{x}^{(k)}_i$. However,~\eqref{eq:local-iterates-inexact} with distinct SSs do not admit at least a weaker form of the consistency, called \emph{synchrony} given by
\begin{equation}\label{eq:Disributed-Sinchrony}
  \forall~i,j\in\mathcal{N},\; i\neq j,~\boldsymbol{x}^{(m)}_j=\boldsymbol{x}^{(m)}_i%, \forall~k\in \Z
\end{equation}
where $m$ is an iteration index of practical interest, e.g., the iteration index at the termination. Thus, the main challenge in this research is to establish the convergence properties of \eqref{eq:local-iterates-inexact}, while maintaining the synchrony.~\footnote{Under imperfect conditions, iterates of the form \eqref{eq:local-iterates-inexact} are commonplace in many distributed algorithms such as primal or dual-decomposition methods, among others, see e.g., \cite{Low-Layering-as-Optimization-Decomposition-2007} and references therein.} This challenge is taken up next, where the iterate~\eqref{eq:local-iterates-inexact} is integrated with potential SS coordination to yield an algorithm with guaranteed convergence.

\section{Algorithm Development}
\label{sec:Algorithm-Development}

Let us first focus on establishing the evolutionary characteristics of~\eqref{eq:local-iterates-inexact} to set the stage for our subsequent developments.

\subsection{Evolutionary Characteristics of~\eqref{eq:local-iterates-inexact}}\label{subsec:Uncontrolled-Worst-Case-Divergence}

From \eqref{eq:local-iterates-inexact}, \eqref{grad_est_model}, together with some standard algebraic manipulations as shown in the \emph{Appendix}, it can be shown that, under \textbf{AS}~\ref{Ass:objectives_assump}, \textbf{AS}~\ref{Ass:Non-Ideal-Communication} and for $\gamma\in(0,1/\sum_{j=i}^{N}L_j]$,
\begin{equation}\label{eq:abs_det_var_gap}
    \| \nabla f\big(\boldsymbol{x}_{i}^{(k)}\big) - \boldsymbol{h}^{(k)}_i\| \leq 2\epsilon N \left( k + {1}/{2}\right), \quad i\in\mathcal{N}
\end{equation}
where $\boldsymbol{h}^{(k)}_i\triangleq\sum_{j=1}^{N} \boldsymbol{h}^{(k)}_{ij}$.  The inequality~\eqref{eq:abs_det_var_gap} indicates that, in the worst case, the norm of the difference between $\nabla f\big(\boldsymbol{x}_{i}^{(k)}\big)$ and its local representation $\boldsymbol{h}_{i}^{(k)}$ diverges as $k\rightarrow\infty$. Thus, it is of paramount importance to control such growth for establishing convergence of iterates of the form~\eqref{eq:local-iterates-inexact}.  To this end, it is customary to rely on SS coordination possibly through an {error-free} communication medium. However, {error-free} communications are usually more expensive. Therefore, unlike the commonly considered periodic SS coordination \cite{Xie2020CSER:Reset_f}, we seek to reduce the communication overhead by dynamically choosing the coordination epochs, so as to make it still possible to ensure convergences of the underlying sequences.
\begin{comment}
to enable {\bf some form of synchrony \color{red} This sounds pretty vague}, see~\cite{Yu2019ParallelLearning_f,Zhou2018OnOptimization_f}. {\bf This in turn affects the communication efficiency, which is critical for distributed...., especially for ...  }However, it is worth pointing out that communication efficiency is of paramount importance in distributed algorithms, especially for machine learning problems of high-dimensional decision spaces~\cite{Mahmoudi2020Cost-efficientNetworks_f,Federated-learning-McMahan-2016,sindri_2018,sindri_2020,Nedic2018NetworkOptimization_f}.

 Therefore, unlike the commonly considered periodic SS coordination, we seek to reduce the communication overhead by dynamically {\bf deciding, \color{red} Can we use an alternative word like "choosing, determining,  "} the coordination epochs, so as to make it still possible to ensure convergences of the underlying sequences.
\end{comment}
%
As such, we consider a relative deviation of the gradient of the objective function $f$ and its measurement from the standpoint of $i$th SS, i.e., $e_i^{(k)}\triangleq\| \nabla f\big(\boldsymbol{x}_{i}^{(k)}\big) - \boldsymbol{h}^{(k)}_i\|/\| \nabla f\big(\boldsymbol{x}_{i}^{(k)}\big)\|,~i\in\mathcal{N}$.

Intuitively, when $e_i^{(k)}$ is sufficiently small, the influence of errors $\boldsymbol{\epsilon}_{ij}^{(k)}$ on \eqref{eq:local-iterates-inexact} becomes relatively insignificant. On the other hand, when $e_i^{(k)}$ is sufficiently large, the consequences become more detrimental, and \eqref{eq:local-iterates-inexact} may evolve anomalously. Thus, to circumvent such anomalies, the objective is to start with synchrony [\cf~\eqref{eq:Disributed-Sinchrony}] at $k=0$ and to perform iterate \eqref{eq:local-iterates-inexact} as long as $e_i^{(k)}$ is sufficiently small, for otherwise to trigger SS coordination. As a result, the iterates \eqref{eq:local-iterates-inexact} at each SSs might tend to evolve in a meaningful direction.

Let us next discuss how the preceding concept can be integrated into devise our algorithm. In this respect, the most crucial step is to identify an epoch at which the SS coordination is to be triggered. In other words, each SS needs a {locally} verifiable characterization of the iterates $k$ for which {$e_i^{(k)}$ is sufficiently small}, despite the dependence of $e_i^{(k)}$ on {global} information $\nabla f\big(\cdot\big)$. As such, we rely on the condition %\begin{equation}\label{eq:maintaining-condition}
%     k\leq \frac{r}{2\epsilon N} \|\boldsymbol{h}_i^{(k)}\|-\frac{1}{2} \implies e_i^{(k)}\leq \frac{r}{1-r}  
% \end{equation}  
\begin{equation}\label{eq:maintaining-condition}
    k\leq {r}\|\boldsymbol{h}_i^{(k)}\|/{(2\epsilon N)} -{1}/{2} \implies e_i^{(k)}\leq {r}/{(1-r)}  
\end{equation} 
\addb{where $r\in(0,1)$ is a design parameter, suitably chosen based on the strong convexity constants and the Lipschitz constants of the objective functions}. The condition~\eqref{eq:maintaining-condition} follows from~\eqref{eq:abs_det_var_gap}, together with that $\|\boldsymbol{h}_i^{(k)}\|-\| \nabla f\big(\boldsymbol{x}_{i}^{(k)}\big)\|\leq \| \nabla f\big(\boldsymbol{x}_{i}^{(k)}\big) - \boldsymbol{h}^{(k)}_i\|$. Thus, the SSs perform the iterate \eqref{eq:local-iterates-inexact} independent of each other, as long as, for all $i\in\mathcal{N}$, $k\leq {r\|\boldsymbol{h}_i^{(k)}\|}/({2\epsilon N})-{1}/{2}$, and is referred to as {\texttt{IndComp}}. If $k> {r\|\boldsymbol{h}_i^{(k)}\|}/({2\epsilon N})-{1}/{2}$ for at least one SS, SSs communicate with others to {average} their local copies $\boldsymbol{x}_{i}^{(k)}$, which is referred to as the intermittent synchronization ({\texttt{IntSync}}). {\texttt{IntSync}} is performed through an {error-free} communication system. Having presented the evolutionary characteristics of ~\eqref{eq:local-iterates-inexact}, we are now ready to propose our new algorithm.
% {\bf Having done blh blah, We are now ready to outline the proposed algorithm.}
\subsection{Algorithm and Its Convergence Analysis}\label{susec:Algorithm-and-Its-Convergence-Analysis}
% Our goal is to develop a systematic switching between {\texttt{IndComp}} and {\texttt{IntSync}} steps, in addition to being convergent, can take its place as a mean of reducing the communication cost under the following hypothesis.
% \begin{remark}
% \label{Remark:Expensive-Inexpensive-Relevance}
% From a practical point of view, {\texttt{IndComp}} steps are less expensive since bounded errors are allowed. However, {\texttt{IntSync}} are expensive since errors are not allowed.
% \end{remark}
The two stages {\texttt{IndComp}} and {\texttt{IntSync}} are implemented in an iterative manner to yield the following algorithm.
\begin{breakablealgorithm} 
\caption{Inexact GD with \texttt{IndComp}$-$\texttt{IntSync}}
\begin{algorithmic}[1]
\addb{\Require{$\boldsymbol{x}^{(0)}_j=\boldsymbol{x}^{(0)}_i$ $\forall~i,j\in\mathcal{N}$, \addb{$\epsilon\geq 0$, $r\in(0,\sqrt{\ell}/(\sqrt{L}+\sqrt{\ell}))$,} $s=0,k=0$}}
\Repeat 
\Repeat 
        \State $\forall~i\in\mathcal{N}$, compute {$\boldsymbol{x}_{i}^{(s{+}k{+}1)}$ from \eqref{eq:local-iterates-inexact}, $k\gets k+1$}
   \Until{$\exists~i\in\mathcal{N}, \ k-1> {r\|\boldsymbol{h}_i^{(s+k-1)}\|}/({2\epsilon N})-{1}/{2}$}
   \If{$k \neq 1$}
   \State $\forall~i\in\mathcal{N}$, $\boldsymbol{x}^{(s+k-1)}_i \gets \frac{1}{N}\textstyle \sum_{j=1}^{N}\boldsymbol{x}_{j}^{(s+k-1)}$
   \State $s\gets s+k-1$, $k\gets 0$
   \Else
    \State $\forall~i\in\mathcal{N}$, $\boldsymbol{x}^{(s+k)}_i \gets \frac{1}{N}\textstyle \sum_{j=1}^{N}\boldsymbol{x}_{j}^{(s+k)}$
   \State $s\gets s+k$, $k\gets 0$
   \EndIf
   \Until{\texttt{a stopping criterion true}}
\end{algorithmic}
\label{alg:Inexact-GD-with-IntSync}
\end{breakablealgorithm}
It is worth emphasizing that the indices $s$ and $k$ of the algorithm have an important interpretation. The index $s$ always represents an iteration at which the {synchrony} [see \eqref{eq:Disributed-Sinchrony}] of the local copies of the decision variables is imposed, \cf~step~7, step~10. The inner loop [\cf~steps~2-4] {always starts} with synchrony. Thus, $k$ represents the local iteration index within the inner loop, which is reset every time the {synchrony} is imposed,~\cf~step~7, step~10. Consequently, $s{+}k$ is simply the {global} iteration index of Algorithm~\ref{alg:Inexact-GD-with-IntSync}.
The following Proposition establishes the convergence of Algorithm~\ref{alg:Inexact-GD-with-IntSync}.

\begin{figure*}[t]
\centering
\begin{subfigure}{.49\textwidth}
  \centering
  % include first image
  \includegraphics[height=0.18\textheight]{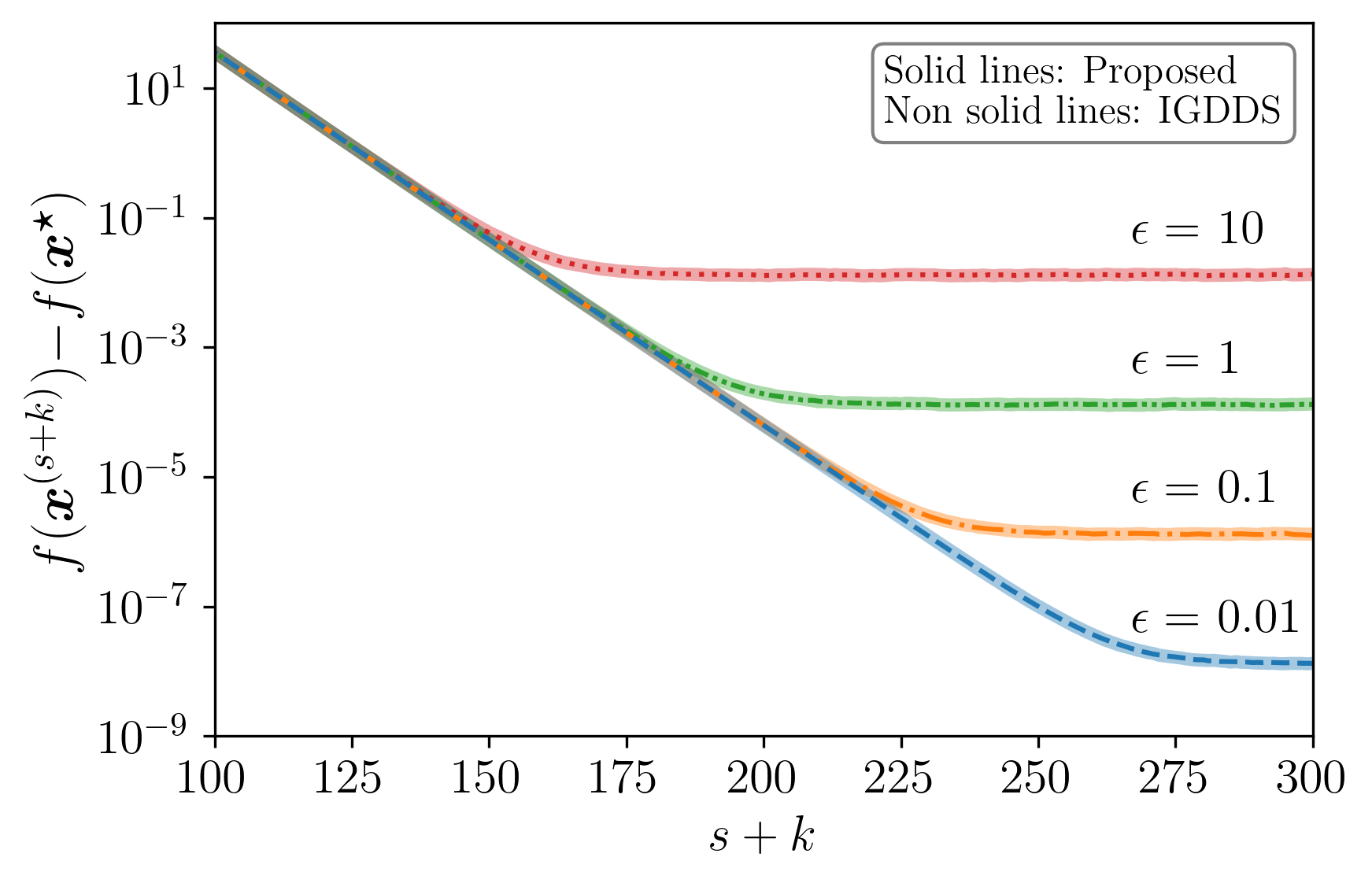} 
  \caption{\footnotesize{Error Vs \texttt{IntSync}  $+$ \texttt{IndComp}}}
  %\caption{}
  \label{fig:simulation_convergence}
\end{subfigure}
\hfill
\begin{subfigure}{.49\textwidth}
  % include second image
\centering
  \includegraphics[height=0.18\textheight]{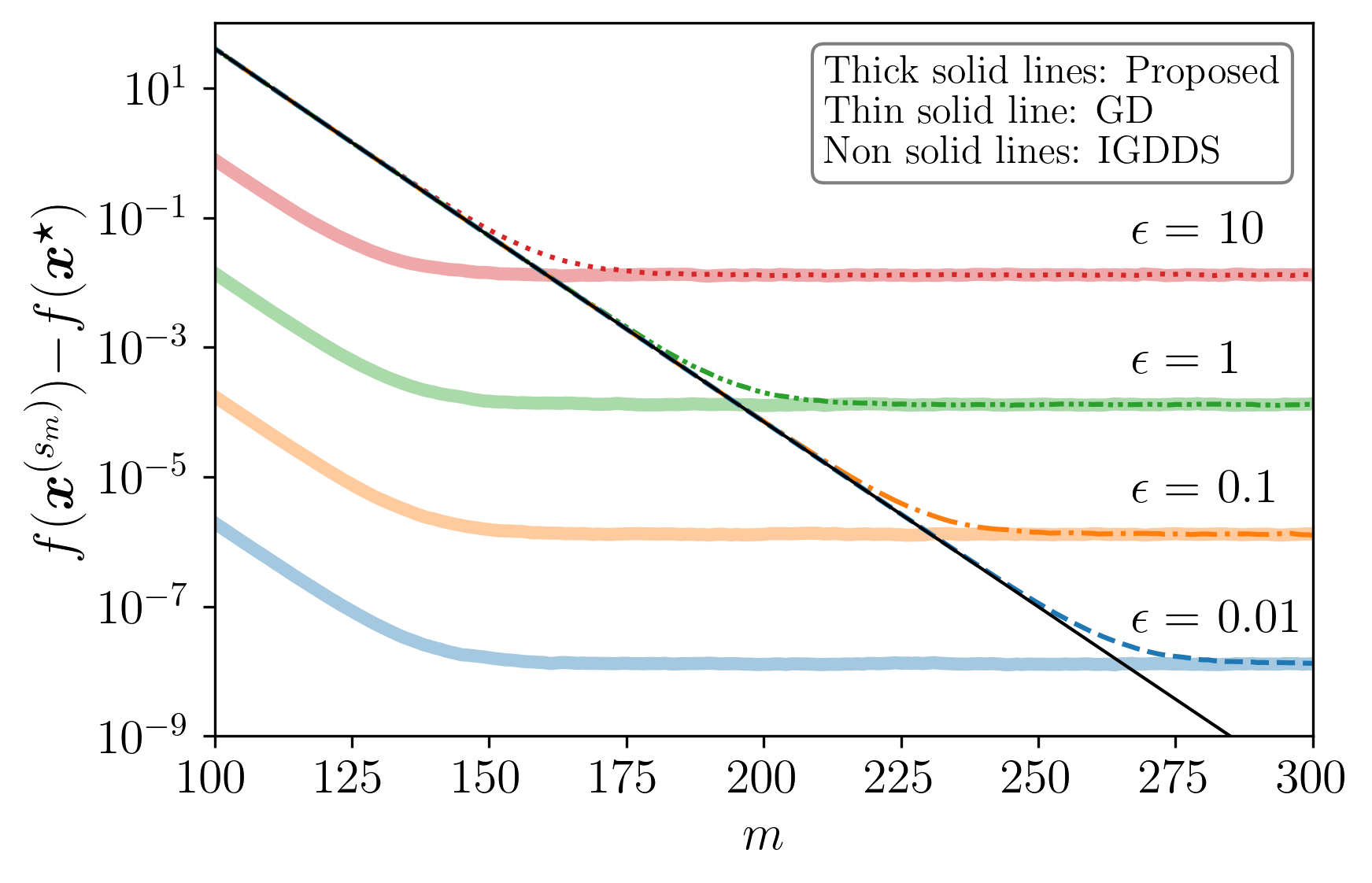}
  \caption{\footnotesize{Error Vs  \texttt{IntSync}}}
  %\caption{}
  \label{fig:simulation_communication_saving}
\end{subfigure}
\caption{Comparison of error for different distortion levels (i.e, $\epsilon$). Results are shown for $\epsilon = 0.01,0.1,1$, and $10$.}

\label{fig:simulation}
\end{figure*}

\begin{prop}
\label{prop:abs_det_grad_bound_c_step}
\addb{Suppose \textbf{AS}~\ref{Ass:objectives_assump}, \textbf{AS}~\ref{Ass:Non-Ideal-Communication} hold. Let $\{\boldsymbol{x}_i^{(k)}\}_{k\in\Z}$, $i\in\mathcal{N}$, be the sequence of local copies of the decision variable generated by Algorithm~\ref{alg:Inexact-GD-with-IntSync}. Then for $\gamma\in(0,1/L]$}% and for %$r\in(0,\sqrt{\ell}/(\sqrt{L}+\sqrt{\ell}))$
   \begin{enumerate}
        \item $\displaystyle \limsup_{k \to \infty} \big(f\big(\boldsymbol{x}^{(k)}_i\big) - f(\boldsymbol{x}^{\star}) \big)\leq {\epsilon^{2} N^{2}}/(2(\ell - L\bar{r}^{2}))$ 
        \item $\displaystyle \limsup_{k \to \infty}\| \nabla f\big(\boldsymbol{x}^{(k)}_i\big)\| \leq \sqrt{{L\epsilon^{2} N^{2}}/(\ell - L \bar{r}^{2})}$
        \item $\displaystyle \limsup_{k \to \infty} \| \boldsymbol{x}^{(k)}_i - \boldsymbol{x}^{\star} \| \leq  \sqrt{{L\epsilon^{2} N^{2}}/(\ell^2 - L \bar{r}^{2}\ell)}$
   \end{enumerate}
where $L=\sum_{j=1}^{N}L_{j}$, $\ell=\min_{j\in\mathcal{N}}\ell_{j}$, $\bar r=r/(1-r)$, and $\boldsymbol{x^\star}{=}\arg\min_{\boldsymbol{x}} f(\boldsymbol x)$. 
\end{prop}

\addb{It is not difficult to see that the Proposition holds even if $\boldsymbol{x}_{i}^{(k)}$ is set as $\boldsymbol{x}^{(k)}=\frac{1}{N}\textstyle \sum_{j=1}^{N}\boldsymbol{x}_{j}^{(k)}$ for all $k\in\Z$.} Note that until the termination of the algorithm [\cf~step~12], the inner loop is in either of the following {states}: 1) it repeats more than once 2) it repeats only once. Thus, the proof of the Proposition is simply based on the characterization of the evolution of the sequence $\{f\big(\boldsymbol{x}_{i}^{(s+k+1)}\big)- f\big(\boldsymbol{x}^\star\big)\}$
%\footnote{\addb{Using Jensen's inequality, one can show that the same results hold for the global average $\boldsymbol{x}^{(k)}$ (i.e: $\boldsymbol{x}^{(k)}=\frac{1}{N}\textstyle \sum_{i=1}^{N}\boldsymbol{x}_{i}^{(k)}$)}).}
when the algorithm is in either of the states. To this end, we shall require the following results, the proofs of which are given in the \emph{Appendix}.
\begin{lemma}
\label{Lemma:First-State}
Let \textbf{AS}~\ref{Ass:objectives_assump}, \textbf{AS}~\ref{Ass:Non-Ideal-Communication} hold, $s\in\Z$ be any iteration index at which synchrony is imposed, $i\in\mathcal{N}$, and $r\in(0,1)$. Moreover, suppose the inner loop of Algorithm~\ref{alg:Inexact-GD-with-IntSync} repeats for iteration indices $\bar k\in\{s,s+1,\ldots,s+\kappa\}$, for some $\kappa\geq 2$. Then for $k\in\{0,1,\ldots,\kappa-1\}$
\begin{equation}
\label{eq:convergence-state-1} 
 \hspace{-1mm}  f\big(\boldsymbol{x}_{i}^{(s+k+1)}\big)- f\big(\boldsymbol{x}^\star\big) \leq q\left(f\big(\boldsymbol{x}_{i}^{(s+k)}\big)-f\big(\boldsymbol{x}^\star\big)\right)
 \end{equation}
 %$L{=}\sum_{j=1}^{N}L_{j}$, 
where $q = (1{+}\gamma L \bar r^2{-}\gamma\ell)$ is a positive constant.
\end{lemma}
% \begin{proof}
% The proof is given in the supplemental material. 
% \end{proof}
Lemma~\ref{Lemma:First-State} characterizes the evolution of the sequence $\{f\big(\boldsymbol{x}_{i}^{(s+k+1)}\big)- f\big(\boldsymbol{x}^\star\big)\}$ when the algorithm is in states~$1$. Consequently, the recursive application of \eqref{eq:convergence-state-1}, together with the Jensen's inequality yields
 \begin{equation}
\label{eq:State-1-Recursive-Application} 
 \hspace{-1mm}  f\big(\boldsymbol{x}_{i}^{(s+\kappa)}\big)- f\big(\boldsymbol{x}^\star\big) \leq q^{\kappa}\left(f\big(\boldsymbol{x}_{i}^{(s)}\big) - f\big(\boldsymbol{x}^\star\big)\right).
 \end{equation} % replace eq label where Corolary:State-1-with-Jensens referred
 %The inequality \eqref{eq:State-1-Recursive-Application} holds even if $\boldsymbol{x}_{i}^{(s+\kappa)}$ from the inner loop of Algorithm~\ref{alg:Inexact-GD-with-IntSync} is set as $\boldsymbol{x}_{i}^{(s+\kappa)}=\frac{1}{N}\textstyle \sum_{j=1}^{N}\boldsymbol{x}_{j}^{(s+\kappa)}$.
 % As a corollary of Lemma~\ref{Lemma:First-State}, we have the following result.
% \begin{corr}
% \label{Corolary:State-1-with-Jensens}
% Let the suppositions of Lemma~\ref{Lemma:First-State} hold. Then recursively application of \eqref{eq:convergence-state-1} for $k\in\{0,1,\ldots,\kappa-1\}$ yields
%  \begin{equation}
% \label{eq:State-1-Recursive-Application} 
%  \hspace{-1mm}  f\big(\boldsymbol{x}_{i}^{(s+\kappa)}\big)- f\big(\boldsymbol{x}^\star\big) \leq q^{\kappa}\left(f\big(\boldsymbol{x}_{i}^{(s)}\big) - f\big(\boldsymbol{x}^\star\big)\right).
%  \end{equation}
%  The inequality \eqref{eq:State-1-Recursive-Application} holds even if $\boldsymbol{x}_{i}^{(s+\kappa)}$ from the inner loop of Algorithm~\ref{alg:Inexact-GD-with-IntSync} is set as $\boldsymbol{x}_{i}^{(s+\kappa)}=(1/N)\textstyle \sum_{j=1}^{N}\boldsymbol{x}_{j}^{(s+\kappa)}$.
% \end{corr}
% \begin{proof}
% The proof is given in the supplemental material. 
% \end{proof}
The evolution of the sequence $\{f\big(\boldsymbol{x}_{i}^{(s+k+1)}\big)- f\big(\boldsymbol{x}^\star\big)\}$ when the algorithm is in state~$2$ is established by the following result.
\begin{lemma}
\label{Lemma:Second-State}
Let \textbf{AS}~\ref{Ass:objectives_assump}, \textbf{AS}~\ref{Ass:Non-Ideal-Communication} hold, $s\in\Z$ be any iteration index at which synchrony is imposed, $r\in(0,1)$, and $i\in\mathcal{N}$. Moreover, suppose the inner loop of Algorithm~\ref{alg:Inexact-GD-with-IntSync} repeats only once, where the iteration index is $s$. Then %for $\gamma\in(0,1/L]$
\begin{align}
\label{eq:convergence-state-2} 
    f\big(\boldsymbol{x}_{i}^{(s+1)}\big) {-} f\big(\boldsymbol{x}^\star \big)  &\leq \big(1 {-} \gamma\ell \big)\left(f\big(\boldsymbol{x}_{i}^{(s)} \big)  - f\big(\boldsymbol{x}^\star \big)\right) \allowdisplaybreaks  
     + \frac{\gamma \epsilon^{2} N^{2}}{2}.
\end{align}
The inequality \eqref{eq:convergence-state-2} holds even if $\boldsymbol{x}_{i}^{(s+1)}$ from the inner loop of Algorithm~\ref{alg:Inexact-GD-with-IntSync} is set as $\boldsymbol{x}_{i}^{(s+1)}=\frac{1}{N}\textstyle \sum_{j=1}^{N}\boldsymbol{x}_{j}^{(s+1)}$.
\end{lemma}
% \begin{proof}
% The proof is given in the supplemental material. 
% \end{proof}
Finally, the following Lemma asserts that the algorithm necessarily switches to state $2$ from state~$1$.
\begin{lemma}
\label{lemma:state_1_to_state_2_transition}
    Let \textbf{AS}~\ref{Ass:objectives_assump}, \textbf{AS}~\ref{Ass:Non-Ideal-Communication} hold. Moreover, suppose $\forall~i\in\mathcal{N}$, $r\|\boldsymbol{h}_i^{(0)}\| \geq \epsilon N$, and thus, the algorithm starts at state~$1$, where $r\in(0,{\sqrt{\ell}}/({\sqrt{L}{+}\sqrt{\ell}})$. Then $\exists$ $\bar s,\bar k \in \Z$ such that Algorithm~\ref{alg:Inexact-GD-with-IntSync} switches to state~$2$ from state~$1$, where $\bar s$ is an iteration index at which the {synchrony} is imposed and $\bar k$ is a local iteration index within the inner loop.
\end{lemma}
% \begin{proof}
% The proof is given in the supplemental material. 
% \end{proof}
Having armed with the above results, we are now ready to give the proof of Proposition~\ref{prop:abs_det_grad_bound_c_step}.
%[\cf~Lemma~\ref{Lemma:First-State} and Lemma~\ref{Lemma:Second-State} is the Appendix]. 

\begin{proof}[Proof of Proposition~\ref{prop:abs_det_grad_bound_c_step}]

%When the algorithm is in \emph{state~1} and \emph{state~2}, the evolution of the sequence $\{f\big(\boldsymbol{x}_{i}^{(k+1)}\big)- f\big(\boldsymbol{x}^\star\big)\}$ is dictated by Lemma

% The following lemma characterizes an upperbound on $\{f\big(\boldsymbol{x}_{i}^{(k+1)}\big)- f\big(\boldsymbol{x}^\star\big)\}$ when the algorithm is in \emph{state~2}.

% Finally, we show that the algorithm necessarily switches to state 2 from state~1.

From Lemma~\ref{Lemma:First-State} and \eqref{eq:State-1-Recursive-Application}, for any consecutive sequence of state~$1$, starting at some global iteration index $n\in\Z$ and ending at $n+k\in\Z$, we have 
\begin{comment}
{\color{red} What is $u^{k}$? Do you mean $u$ to the power $k$? If not, there exists a notational inconsistency here. We do not need (5), (6), (7), and (8) to express the very essence of your argument here.} 
{\color{blue} $u^{k}$ means $u$ to power $k$}
\end{comment}
\begin{align}
\label{eq:Proposition-state1-iterates} 
 \hspace{-1mm}  f\big(&\boldsymbol{x}_{i}^{(n+k)}\big)- f\big(\boldsymbol{x}^\star\big)  \leq q^{k}\left(f\big(\boldsymbol{x}_{i}^{(n)}\big) - f\big(\boldsymbol{x}^\star\big)\right)\\ \label{eq:Proposition-state1-iterates-2} 
 & \leq q^{k}\left(f\big(\boldsymbol{x}_{i}^{(n)}\big) - f\big(\boldsymbol{x}^\star\big)\right) + \frac{\gamma \epsilon^{2} N^{2}}{2}\textstyle \sum_{j=0}^{k-1}q^j.
 \end{align}
Similarly, recursively applying \eqref{eq:convergence-state-2} in Lemma~\ref{Lemma:Second-State} for any consecutive sequence of state~$2$, starting at some global iteration index $n\in\Z$ and ending at $n+k\in\Z$, together with that $1-\gamma\ell\leq q$, we again have an equivalent form of \eqref{eq:Proposition-state1-iterates-2}.
% have
% \begin{align} \nonumber
%  \hspace{-1mm}  f\big(&\boldsymbol{x}_{i}^{(n+k)}\big)- f\big(\boldsymbol{x}^\star\big) \\ \label{eq:Proposition-state2-iterates-2} 
%  & \leq u^{k}\left(f\big(\boldsymbol{x}_{i}^{(n)}\big) - f\big(\boldsymbol{x}^\star\big)\right) + \frac{\gamma \epsilon^{2} N^{2}}{2}\textstyle \sum_{j=0}^{k-1}u^j.
%  \end{align}
Moreover, the algorithm necessarily switches to state~$2$ from state~$1$, \cf~Lemma~\ref{lemma:state_1_to_state_2_transition}. Thus, from \eqref{eq:Proposition-state1-iterates-2}, $\forall$ $k\in\Z$, we have
%and thus for all $k\in\Z$ and \eqref{eq:Proposition-state2-iterates-2}, 
\begin{align} \nonumber
 \hspace{-1mm}  f\big(&\boldsymbol{x}_{i}^{(k)}\big){-} f\big(\boldsymbol{x}^\star\big)  \leq q^{k}\left(f\big(\boldsymbol{x}_{i}^{(0)}\big) {-} f\big(\boldsymbol{x}^\star\big)\right) {+} \frac{\gamma \epsilon^{2} N^{2}}{2}\textstyle \sum_{j=0}^{k-1}q^j.
 \end{align}
Noting that $q<1$, we take the limit as $k\to \infty$ to yield Part~1. Part~2 follows from Part~1 and \cite[eq.~10, \S~1.4]{poliak}. Finally, Part~3 follows from Part~1 and \cite[eq.~35, \S~1.1]{poliak}.
\end{proof}

\begin{comment}
\begin{figure}
  \centering
  % include first image
  \includegraphics[width=0.45\textwidth]{Figures/Averaged Plots/python/Comparison of convergence - function value.jpg} 
  \caption{\footnotesize{Error Vs \texttt{IntSync}  $+$ \texttt{IndComp}}}
  \label{fig:simulation_convergence}
\end{figure}
\begin{figure}
  % include second image
  \includegraphics[width=0.45\textwidth]{Figures/Averaged Plots/python/Comparison of convergence - at intsync steps.jpg}
  \caption{\footnotesize{Error Vs  \texttt{IntSync}}}
  \label{fig:simulation_communication_saving}
\end{figure}
\end{comment}

\section{Numerical Results}
\label{sec:num_results}

Let us first verify the convergence results of Proposition~\ref{prop:abs_det_grad_bound_c_step}. \addb{To this end, we consider problem~\eqref{eq:main-problem} with quadratic $f_i$s, i.e., $f_i(\boldsymbol{x})= \boldsymbol{x}\tran \boldsymbol{B}\tran_i\boldsymbol{B}_i\boldsymbol{x}+\boldsymbol{c}_i\tran \boldsymbol{x}$, where $\boldsymbol{B}\tran_i\boldsymbol{B}_i\in \mathbb{S}^{n}_{++}$, $\boldsymbol{c}_i\in \mathbb{R}^n$, and $\mathbb{S}^{n}_{++}$ is the {positive definite cone}.
% \begin{equation} \label{eq:Numerical-section-local-problems}
%     f_i(\boldsymbol{x})= \boldsymbol{x}\tran \boldsymbol{A}_i \boldsymbol{x}+\boldsymbol{c}_i\tran \boldsymbol{x}, \quad \boldsymbol{A}_i\in \mathbb{S}^{n}_{++}, \ \boldsymbol{c}_i\in \mathbb{R}^n,
% \end{equation}
% where $\mathbb{S}^{n}_{++}$ is the \emph{positive definite cone}. 
The entries of $\boldsymbol{B}_i$ and $\boldsymbol{c}_i$ are generated from a normal distribution. Note that $\ell_i$ and $L_i$ are determined by $\boldsymbol{B}_i$, \cf~\textbf{AS}~\ref{Ass:objectives_assump}. We let $N=4$, $n=10$, $\gamma= 1/(2L)$, and $r = 0.03$.} Only the results related to Proposition~\ref{prop:abs_det_grad_bound_c_step}-(1)  is presented, since those related to Proposition~\ref{prop:abs_det_grad_bound_c_step}-(2) and (3) behave similarly. 

For comparison, we consider two algorithms. The first one is the classic GD, i.e.,~Algorithm~\ref{alg:Inexact-GD-with-IntSync} with $\epsilon = 0$ and $r=0$. We also consider another algorithm which we refer to as inexact-GD with distributed synchrony (IGDDS), i.e.,~Algorithm~\ref{alg:Inexact-GD-with-IntSync} with $r=0$ and $\forall~i\in\mathcal{N},~\boldsymbol\epsilon_{ij}^{(k)}=\boldsymbol\epsilon_{j}^{(k)}$ [\cf~\eqref{grad_est_model}]. In this respect, the synchrony~\eqref{eq:Disributed-Sinchrony} holds for all $k\in\Z$ and we have $ \limsup_{k \to \infty} \big(f\big(\boldsymbol{x}^{(k)}_i\big) - f(\boldsymbol{x}^{\star}) \big)\leq {\epsilon^{2} N^{2}}/(2\ell)$~\cite[\S~4]{Ajalloeian2020OnGradients_f,poliak}. 

Figure~\ref{fig:simulation_convergence} shows the error $f\big(\boldsymbol{x}^{(s+k)}\big){-}f(\boldsymbol x^\star)$ vs global iteration index $s+k$ for different $\epsilon$, \cf~solid lines. \addb{Results are averaged over $1000$ initializations $\boldsymbol{x}^{(0)}$, whose entries are normally distributed.} Plots agree with Proposition~\ref{prop:abs_det_grad_bound_c_step}-(1), i.e., the smaller the $\epsilon$, the smaller the error of the optimality. Results with IGDDS are given in {non-solid} lines. Convergence rates and the suboptimality obtained by Algorithm~\ref{alg:Inexact-GD-with-IntSync} and IGDDS seem almost identical. This is expected because the convergence rate of Algorithm~\ref{alg:Inexact-GD-with-IntSync}, i.e., $(1-\gamma L\bar r^2-\gamma\ell)$ and that of IGDDS, i.e.,~$(1-\gamma\ell)$ are almost identical when $\gamma L \bar r^2\ll 1-\gamma\ell$. This condition is always realizable in practice, e.g., we have $\gamma L \bar r^2=0.0005$ and $1-\gamma\ell=0.9986$ in our simulation. A similar comparison holds for the suboptimality as well. Thus, results suggest that Algorithm~\ref{alg:Inexact-GD-with-IntSync} yields almost identical results to that of more constrained IGDDS. 

Since IGDDS is technically equivalent to Algorithm~\ref{alg:Inexact-GD-with-IntSync} with $r=0$, error-free communication is needed in every iteration to yield synchrony~\eqref{eq:Disributed-Sinchrony}. However, Algorithm~\ref{alg:Inexact-GD-with-IntSync} does not require synchrony in every iteration. Therefore, for a {fair} comparison of Algorithm~\ref{alg:Inexact-GD-with-IntSync} and IGDDS in terms of communication overhead, it is instructive to plot the error versus the number of \texttt{IntSync} steps $m$, where $s_m$, ${m\in\Z}$ is the iteration index of $s$ within Algorithm~\ref{alg:Inexact-GD-with-IntSync} at which the $m$th-synchrony is imposed.

Figure~\ref{fig:simulation_communication_saving} shows the error $f\big(\boldsymbol{x}^{(s_m)}\big)-f(\boldsymbol x^\star)$ vs $m$ with Algorithm~\ref{alg:Inexact-GD-with-IntSync}, see thick solid lines. Results related to IGDDS are also plotted, see the {non-solid} lines. Clearly, there is a shift of the plots with IGDDS towards the right relative to the plots with Algorithm~\ref{alg:Inexact-GD-with-IntSync}. Therefore, for all considered $\epsilon$ values, the number of \texttt{IntSync} steps $m$ required to obtain a specified error with Algorithm~\ref{alg:Inexact-GD-with-IntSync} is smaller than with IGDDS. Moreover, if the number of \texttt{IntSync} steps $m$ is fixed, the error with Algorithm~\ref{alg:Inexact-GD-with-IntSync} can be on the order of magnitude smaller than with IGDDS. This is useful in practice, because the cost of the error-free communication required for \texttt{IntSync} can be reduced with Algorithm~\ref{alg:Inexact-GD-with-IntSync} than with IGDDS. %[\cf~Remark~\ref{Remark:Expensive-Inexpensive-Relevance}]. 
The benefits become greater as $\epsilon$ decreases. Finally, we plot results due to GD, see the thin solid line in Fig.~\ref{fig:simulation_communication_saving}. Results show that still the Algorithm~\ref{alg:Inexact-GD-with-IntSync} can benefit from less expensive \texttt{IndComp} steps. %[\cf Remark~\ref{Remark:Expensive-Inexpensive-Relevance}]. 
For example, in $150$ \texttt{IntSync} steps, Algorithm~\ref{alg:Inexact-GD-with-IntSync} manages to yield an error significantly less than that from GD despite the value of $\epsilon$. Clearly, GD outperforms Algorithm~\ref{alg:Inexact-GD-with-IntSync} if $m$ is sufficiently large, since there are no inexactnesses. Thus, the results suggest if there is a choice for less expensive communication for \texttt{IndComp}, or a choice for allowing some inexactnesses, one can operate Algorithm~\ref{alg:Inexact-GD-with-IntSync} in a way there is a trade-off between the error and \texttt{IntSync} steps ($m$).

\section{Conclusion} \label{sec:conclusion}
A gradient-like algorithm with guaranteed convergence has been developed to minimize a sum of peer objective functions through an interconnection network with multi-peer broadcast and multi-peer accumulation capabilities. Peer coordination can usually admit communications with bounded errors, however with some infrequent error-free synchronization epochs, which are dynamically triggered. Our algorithm can be attractive in many distributed applications, under inexact communication settings, such as decomposition with dual-subgradient methods and distributed learning systems with in-network computing capabilities, among others.

\clearpage
\appendix
\subsection{Derivation of \eqref{eq:abs_det_var_gap}}
\label{s_sec:derrivation_abs_det_var_gap}
Suppose \textbf{AS}~\ref{Ass:objectives_assump} and \textbf{AS}~\ref{Ass:Non-Ideal-Communication} hold. Moreover, let $\gamma\in(0,1/\sum_{j=i}^{N}L_j]$. Then by recursively applying \eqref{eq:local-iterates-inexact} followed by the use of the triangular inequality gives
\begin{equation}
\label{eq:x_gap_s}
    \|\boldsymbol{x}_{i}^{(k)} - \boldsymbol{x}_{j}^{(k)}\| \leq 2\epsilon N\gamma k, \;\; \forall i,j\in\mathcal{N}.
\end{equation}
 From~\eqref{grad_est_model} and the gradient Lipshitz continuity of $f_i$s, it follows that
\begin{equation}\label{Divergence-grad}
    \| \nabla f_j\big(\boldsymbol{x}_{i}^{(k)}\big) - \boldsymbol{h}_{ij}^{(k)}\|\leq L_j\|\boldsymbol{x}_{i}^{(k)}-\boldsymbol{x}_{j}^{(k)}\| + \epsilon.
\end{equation}
Finally, \eqref{eq:abs_det_var_gap} follows from \eqref{Divergence-grad} by noting that $\nabla f\big(\boldsymbol{x}^{(k)}_i\big)=\sum_{j=1}^{N} \nabla f_{j}\big(\boldsymbol{x}^{(k)}_j\big)$, the definition of $\boldsymbol{h}^{(k)}_i$, \eqref{eq:x_gap_s}, and $\gamma\leq \left(1/\sum_{j=i}^{N}L_j\right)$.

\subsection{Proof of Lemma~\ref{Lemma:First-State}}
\label{s_sec:proof_lemma_first_state}
Without loss of generality we may assume $s=0$. Now, one can bound $f(\boldsymbol{x}_i^{(k+1)})$ as follows:
\begin{align}  \nonumber % \label{descent_lemma} 
%\begin{split}
  \hspace{-3mm}  f\big(\boldsymbol{x}_{i}^{(k+1)}\big)& \leq f\big(\boldsymbol{x}_{i}^{(k)}\big) {+} \nabla f\big(\boldsymbol{x}_{i}^{(k)}\big)\tran\big(\boldsymbol{x}_{i}^{(k+1)} - \boldsymbol{x}_{i}^{(k)}\big) \allowdisplaybreaks \\  
    & \hspace{30mm} ({L}/{2})\|\boldsymbol{x}_{i}^{(k+1)} - \boldsymbol{x}_{i}^{(k)}\|^{2}  \label{eq:convergence-1}  \allowdisplaybreaks \\
    & \leq f\big(\boldsymbol{x}_{i}^{(k)}\big) {-}\gamma \nabla f\big(\boldsymbol{x}_{i}^{(k)}\big)\tran\boldsymbol{h}^{(k)}_i{+}({\gamma}/{2})\|\boldsymbol{h}^{(k)}_i\|^2 \label{eq:convergence-2} \allowdisplaybreaks \\
    & = f\big(\boldsymbol{x}_{i}^{(k)}\big) {-}({\gamma}/{2}) \|\nabla f\big(\boldsymbol{x}_{i}^{(k)}\big)\|^2 \allowdisplaybreaks \nonumber \\ 
    &  \hspace{28mm} ({\gamma}/{2}) \|\nabla f\big(\boldsymbol{x}_{i}^{(k)}\big)-\boldsymbol{h}^{(k)}_i\|^2 \label{eq:convergence-3}  \allowdisplaybreaks \\
    & \leq  f\big(\boldsymbol{x}_{i}^{(k)}\big) {+}\left(\frac{\gamma{\bar r}^2}{2}-\frac{\gamma}{2}\right) \|\nabla f\big(\boldsymbol{x}_{i}^{(k)}\big)\|^2 \label{eq:convergence-4}  \allowdisplaybreaks \\
     & \leq  f\big(\boldsymbol{x}_{i}^{(k)}\big) +(\gamma L{\bar r}^2{-}{\gamma}{\ell}) \left(f\big(\boldsymbol{x}_{i}^{(k)}\big)-f\big(\boldsymbol{x}^\star\big)\right) \label{eq:convergence-5}
   % \end{split}
\end{align}
where $\bar r=r/(1-r)$. Here \eqref{eq:convergence-1} follows from the \emph{descent lemma}~\cite[Lemma 5.7]{beck2017first}, \eqref{eq:convergence-2} follows from \eqref{eq:local-iterates-inexact} and noting that $\gamma L\leq 1$, \eqref{eq:convergence-3} follows from simple algebraic identities, \eqref{eq:convergence-4} follows from \eqref{eq:maintaining-condition}, \eqref{eq:convergence-5} follows from \cite[Lemma~3, \S~1.4]{poliak} and \cite[eq.~10, \S~1.4]{poliak} for bounding $-\|\nabla f\big(\boldsymbol{x}_{i}^{(k)}\big)\|^2$ and $\|\nabla f\big(\boldsymbol{x}_{i}^{(k)}\big)\|^2$, respectively. Now, subtracting $f\big(\boldsymbol{x}^\star\big)$ from the both sides of \eqref{eq:convergence-5} yields the final result.

\subsection{Proof of Lemma~\ref{Lemma:Second-State}}
\label{s_sec:proof_lemma_second_state}
To begin with, let us bound $f(\boldsymbol{x}_i^{(s+1)})$ as follows:
\begin{align}  \nonumber % \label{descent_lemma} 
%\begin{split}
  \hspace{-3mm}  f\big(\boldsymbol{x}_{i}^{(s+1)}\big)& \leq f\big(\boldsymbol{x}_{i}^{(s)}\big) {-}({\gamma}/{2}) \|\nabla f\big(\boldsymbol{x}_{i}^{(s)}\big)\|^2 \allowdisplaybreaks \nonumber \\ 
    &  \hspace{20mm} + ({\gamma}/{2}) \|\nabla f\big(\boldsymbol{x}_{i}^{(s)}\big)-\boldsymbol{h}^{(s)}_i\|^2  \label{eq:state_2_convergence-1} \allowdisplaybreaks \\
    & \leq  f\big(\boldsymbol{x}_{i}^{(s)}\big) - \gamma \ell \left(f\big(\boldsymbol{x}_{i}^{(s)}\big)-f\big(\boldsymbol{x}^\star\big)\right)  \allowdisplaybreaks 
     {+} \big(\gamma \epsilon^{2} N^{2}\big)/2 
    \label{eq:state_2_convergence-2} 
   % \end{split}
\end{align}
where \eqref{eq:state_2_convergence-1} is similar to \eqref{eq:convergence-3} of the preceding lemma. \eqref{eq:state_2_convergence-2} follows from \cite[Lemma~3, \S~1.4]{poliak} for bounding $-\|\nabla f\big(\boldsymbol{x}_{i}^{(s)}\big)\|^2$ and from that $\|\nabla f\big(\boldsymbol{x}_{i}^{(s)}\big)-\boldsymbol{h}^{(s)}_i\|^2\leq \epsilon N$, since the inner loop always starts from synchrony, \cf~\eqref{eq:abs_det_var_gap}. Subtracting $f\big(\boldsymbol{x}^\star\big)$ from both sides yields \eqref{eq:convergence-state-2}. The latter part of the lemma is immediate from the Jensen's inequality.

\subsection{Proof of Lemma~\ref{lemma:state_1_to_state_2_transition}}
\label{s_sec:proof_lemma_state_1_to_state_2_transition}
Suppose the algorithm remains in state~$1$~\footnote{More generally, the algorithm can be in a consecutive sequence of inner loops that are of state~$1$.} without switching to state~$2$. For clarity, let $s\in\Z$ and $k\in\Z$ denote arbitrary iteration indices at which the \emph{synchrony} is imposed and corresponding local iteration index within the inner loop, respectively. Thus, from Lemma~\ref{Lemma:First-State} and \eqref{eq:State-1-Recursive-Application}, we have
 \begin{equation}
\label{eq:state_1_to_state_2_switch-1} 
 \hspace{-1mm}  f\big(\boldsymbol{x}_{i}^{(s+k)}\big) - f\big(\boldsymbol{x}^\star\big) \leq q^{s+k}\left(f\big(\boldsymbol{x}_{i}^{(0)}\big) - f\big(\boldsymbol{x}^\star\big)\right).
 \end{equation}
Consequently, bounding $f\big(\boldsymbol{x}_{i}^{(s+k)}\big)-f\big(\boldsymbol{x}^\star\big)$ and $f\big(\boldsymbol{x}_{i}^{(0)}\big)- f\big(\boldsymbol{x}^\star\big)$ using  \cite[Lemma~3, \S~1.4]{poliak} and \cite[eq.~10 \S~1.4]{poliak} respectively, we have 
 \begin{equation}
\label{eq:state_1_to_state_2_switch-2}  
 \hspace{-1mm}  \|\nabla f\big(\boldsymbol{x}_{i}^{(s+k)}\big)\|^{2} \leq ({L}/{\ell})q^{s+k}\|\nabla f\big(\boldsymbol{x}_{i}^{(0)}\big)\|^{2}.
 \end{equation}
Moreover, for $r\in(0,{\sqrt{\ell}}/({\sqrt{L}{+}\sqrt{\ell}}))$, we have $ q\in(0,1)$. Thus, $\exists~ s+k \in \mathbb{Z}^{+}$ such that guarantees \footnote{To be precise $ s+k \geq \left\lceil \frac{\ln (\epsilon^2N^2L\|\nabla f\big(\boldsymbol{x}_{i}^{(0)}\big)\|^{2})-\ln(\bar r^2 \ell)}{\ln q }\right\rceil$, where $ \lceil . \rceil$ is the ceiling function.}
\begin{equation}
 \label{eq:state_1_to_state_2_switch-4} 
     \|\nabla f\big(\boldsymbol{x}_{i}^{(s+k)}\big)\| < {\epsilon N}/{\bar{r}}.
 \end{equation}
It holds that
\begin{align}  \label{eq:state_1_to_state_2_switch-5} 
 \hspace{-3mm}   \|\boldsymbol{h}_i^{(s+k)}\| & \leq \|\nabla f\big(\boldsymbol{x}_{i}^{(s+k)}\big)\| {+} \| \nabla f\big(\boldsymbol{x}_{i}^{(s+k)}\big) - \boldsymbol{h}^{(s+k)}_i\| \\  \label{eq:state_1_to_state_2_switch-6} 
 & < {\epsilon N}/{\bar{r}}+ 2 \epsilon N \big(k + 1/2 \big)\\  \label{eq:state_1_to_state_2_switch-7}  
 &< 2{\epsilon N}(k+1/2)\big(1/\bar r+1\big)={2{\epsilon N}(k{+}1/2)}/{r}
\end{align}
where \eqref{eq:state_1_to_state_2_switch-5} follows from triangular inequality, \eqref{eq:state_1_to_state_2_switch-6} follows from \eqref{eq:state_1_to_state_2_switch-4} and \eqref{eq:abs_det_var_gap}.

The inequality \eqref{eq:state_1_to_state_2_switch-7} is the inner loop exit criterion [\cf~step~4] which transfers the control of the algorithm to \texttt{IntSync} at steps~6-7 of the algorithm. From \eqref{eq:State-1-Recursive-Application} it follows that the inequality \eqref{eq:state_1_to_state_2_switch-1} holds even after the synchrony at \texttt{IntSync}. Thus, by following arguments identical to that of \eqref{eq:state_1_to_state_2_switch-2} - \eqref{eq:state_1_to_state_2_switch-7}, we conclude that the control of the algorithm is next transferred to \texttt{IntSync} at steps~9-10. That is, the previous inner loop has been repeated only once, which is a contradiction. Therefore, the algorithm must switch to state 2.

\subsection{Analysis with a General Peer-to-Peer Setting}\label{App:Generalized-Graph}
In \S~\ref{sec:prob_formulation} and \S~\ref{sec:Algorithm-Development}, we focused on a network that can be modeled using a fully connected graph. However, the mathematical derivations can be extended to a more generalized peer-to-peer network that is modeled using a connected graph. Therefore, communication need not be coordinated by a central controller like in a federated setting. In the sequel, the main points of the derivations and related results are discussed.

Let us consider an arbitrary graph $\mathcal{G}(\mathcal{N}, \mathcal{L})$, where $\mathcal{N} = \{1,2,\ldots,N\}$ represents the set of subsystems (SSs) of problem~\eqref{eq:main-problem}. Moreover, $\mathcal{L}$ represents a set of edges between SSs, where an edge is given by a pair $(i,j)$, $i,j\in\mathcal{N}$. The graph is considered to be undirected. In other words, $(i,j)\in\mathcal{L}\iff(j,i)\in\mathcal{L}$. Communication from SS~$j$ to $i$ is allowed if and only if there is a link $(i,j)$ between the two nodes. We denote by $\mathcal{N}_{i} = \left\{j \ | \ (i,j) \in \mathcal{L} \right\}$, the set of neighbours of the SS $i$. Furthermore, we denote by $\mathcal{T}(\mathcal{N}, \Bar{\mathcal{L}})$, a spanning tree of the graph $\mathcal{G}$ where $\bar{\mathcal{L}}$ denotes the set of edges in $\mathcal{T}$. We also define $\mathcal{T}_{i} = \left\{j \ | \ (i,j) \in \Bar{\mathcal{L}} \right\}$, the neighbors of the SS $i$ in the spanning tree. 

Now, we note that the gradient measurement model in \eqref{grad_est_model} is going to be modified as follows in the general setting:
\begin{equation}
\boldsymbol{h}_{ij}^{(k)} = \left\{ 
                                \begin{array}{cl}
                                   \nabla f_{j} \big(\boldsymbol{x}_j^{(k)} \big) + \boldsymbol{\epsilon}_{ij}^{(k)}  &  \mbox{if $(i,j) \in \mathcal{L}$} \\ 
                                   \boldsymbol{0}  &  \mbox{otherwise}.
                                \end{array}
                            \right.
\end{equation}
Consequently, it is immediate that
\begin{equation}
\label{eq:measurement-error}
    \left\| \boldsymbol{h}_{ij}^{(k)} -  \nabla f_{j} \big(\boldsymbol{x}_j^{(k)} \big) \right\| = \left\{ 
                                \begin{array}{cl}
                                   \left\| \boldsymbol{\epsilon}_{ij}^{(k)} \right\| &  \mbox{if $(i,j) \in \mathcal{L}$} \\ 
                                   \left\| \nabla f_{j} \big(\boldsymbol{x}_j^{(k)} \big) \right\|  &  \mbox{otherwise}.
                                \end{array}
                            \right.
\end{equation}
%which is the error of the gradient measurement received to SS $i$ from SS $j$. 
It is worth highlighting that the generalized setting requires an additional assumption unlike the fully connected setting considered in \S~\ref{sec:prob_formulation} and \S~\ref{sec:Algorithm-Development}, which we will outline next.
\begin{assumption}
\label{Ass:gradient_boundness}
    The gradients $\nabla f_{i}$s of the objective functions $f_{i}$, $i\in\mathcal{N}$ are bounded, i.e., $\left\|\nabla f_{i}(\boldsymbol{x})\right\| \leq \zeta$ for some $\zeta > 0$ for all $\boldsymbol{x}, i$.
\end{assumption}
Hence, from \textbf{AS}~\ref{Ass:Non-Ideal-Communication} and \textbf{AS}~\ref{Ass:gradient_boundness}, together with \eqref{eq:measurement-error}, it is easily verified that
\begin{equation}
    \left\| \boldsymbol{h}_{ij}^{(k)} -  \nabla f_{j} \big(\boldsymbol{x}_j^{(k)} \big) \right\| \leq \tau 
\end{equation}
where $\tau = \max \big({\epsilon, \zeta}\big)$. Thus, the gradient measurement $\boldsymbol{h}_{ij}^{(k)}$ obeys the following remark:
\begin{remark} 
\label{Rem:Non-Ideal-Communication}
$\forall$ $i,j\in\mathcal{N}$, s.t. $i\neq j$, gradient measurement $\boldsymbol{h}_{ij}^{(k)}\in\R^n$ received by $i$-th SS from $j$-th SS at $k$-th iteration is given by 
\begin{equation}
\label{grad_est_model_new}
    \boldsymbol{h}_{ij}^{(k)} = \nabla f_{j} \big(\boldsymbol{x}_j^{(k)} \big) + \boldsymbol{\tau}_{ij}^{(k)}
\end{equation}
where $\boldsymbol{\tau}_{ij}^{(k)}\in\R^n$ is a error such that {$||\boldsymbol{\tau}_{ij}^{(k)}||\leq \tau$ with $||\cdot||$ denoting the Euclidean norm.}
\end{remark} 
Note that Remark~\ref{Rem:Non-Ideal-Communication} play the role of \textbf{AS}~\ref{Ass:Non-Ideal-Communication}. 

Let us next outline the modified version of Algorithm~\ref{alg:Inexact-GD-with-IntSync}. 

% \begin{breakablealgorithm} 
% \caption{Inexact GD with \texttt{IndComp}$-$\texttt{IntSync}(Modified)}
% \begin{algorithmic}[1]
% \Require{$\boldsymbol{x}^{(0)}_j=\boldsymbol{x}^{(0)}_i$ $\forall~i,j\in\mathcal{N}$, $\epsilon$, $r$, $s=0,k=0$} 
% \Repeat 
% \Repeat 
%         \State $\forall~i\in\mathcal{N}$, compute {$\boldsymbol{x}_{i}^{(s{+}k{+}1)}$ from (3), $k\gets k+1$}
%    \Until{$\exists~i\in\mathcal{N}, \ k-1> {r\|\boldsymbol{h}_i^{(s+k-1)}\|}/({2\epsilon N})-{1}/{2}$}
%    \If{$k \neq 1$}
%    \State \texttt{IntSync}($s+k-1$)
%    \State $s\gets s+k-1$, $k\gets 0$
%    \Else
%     \State \texttt{IntSync}($s+k$)
%    \State $s\gets s+k$, $k\gets 0$
%    \EndIf
%    \Until{\texttt{a stopping criterion true}}
% \end{algorithmic}
% \label{alg:Inexact-GD-with-IntSync-modified}
% \end{breakablealgorithm}
\begin{breakablealgorithm} 
\caption{Inexact GD with \texttt{IndComp}$-$\texttt{IntSync} over a General Graph}
\begin{algorithmic}[1]
\addb{\Require{$\boldsymbol{x}^{(0)}_j=\boldsymbol{x}^{(0)}_i$ $\forall~i,j\in\mathcal{N}$, \addb{$\tau\geq 0$, $r\in(0,\sqrt{\ell}/(\sqrt{L}+\sqrt{\ell}))$,} $s=0,k=0$}}
\Repeat 
\Repeat 
        \State $\forall~i\in\mathcal{N}$, compute {$\boldsymbol{x}_{i}^{(s{+}k{+}1)}$ from (3), $k\gets k+1$}
   \Until{$\exists~i\in\mathcal{N}, \ k-1> {r\|\boldsymbol{h}_i^{(s+k-1)}\|}/({2\tau N})-{1}/{2}$}
   \If{$k \neq 1$}
   \State $\forall~i\in\mathcal{N}$, $\boldsymbol{x}^{(s+k-1)}_i \gets \frac{1}{N}\textstyle \sum_{j=1}^{N}\boldsymbol{x}_{j}^{(s+k-1)}$ \Comment{performed through $\mathcal{T}(\mathcal{G}, \Bar{\mathcal{L}})$.}
   \State $s\gets s+k-1$, $k\gets 0$
   \Else
    \State $\forall~i\in\mathcal{N}$, $\boldsymbol{x}^{(s+k)}_i \gets \frac{1}{N}\textstyle \sum_{j=1}^{N}\boldsymbol{x}_{j}^{(s+k)}$ \Comment{performed through $\mathcal{T}(\mathcal{G}, \Bar{\mathcal{L}})$.}
   \State $s\gets s+k$, $k\gets 0$
   \EndIf
   \Until{\texttt{a stopping criterion true}}
\end{algorithmic}
\label{alg:Inexact-GD-with-IntSync-2}
\end{breakablealgorithm}
Now, one can easily see that an identical result to Proposition~\ref{prop:abs_det_grad_bound_c_step} holds even in the general setting if \textbf{AS}~\ref{Ass:Non-Ideal-Communication} is replaced by Remark~\ref{Rem:Non-Ideal-Communication} above. More specifically, we have the following result:
\begin{prop}
\label{prop:abs_det_grad_bound_c_step_extended}
Suppose \textbf{AS}~1, Remark~\ref{Rem:Non-Ideal-Communication} hold. 
Let $\{\boldsymbol{x}_i^{(k)}\}_{k\in\Z}$, $i\in\mathcal{N}$, be the sequence of local copies of the decision variable generated by Algorithm~\ref{alg:Inexact-GD-with-IntSync-2}. Then for $\gamma\in(0,1/L]$% and for 
   \begin{enumerate}
        \item $\displaystyle \limsup_{k \to \infty} \big(f\big(\boldsymbol{x}^{(k)}_i\big) - f(\boldsymbol{x}^{\star}) \big)\leq {\tau^{2} N^{2}}/(2(\ell - L\bar{r}^{2}))$ 
        \item $\displaystyle \limsup_{k \to \infty}\| \nabla f\big(\boldsymbol{x}^{(k)}_i\big)\| \leq \sqrt{{L\tau^{2} N^{2}}/(\ell - L \bar{r}^{2})}$
        \item $\displaystyle \limsup_{k \to \infty} \| \boldsymbol{x}^{(k)}_i - \boldsymbol{x}^{\star} \| \leq  \sqrt{{L\tau^{2} N^{2}}/(\ell^2 - L \bar{r}^{2}\ell)}$
   \end{enumerate}
where $L=\sum_{j=1}^{N}L_{j}$, $\ell=\min_{j\in\mathcal{N}}\ell_{j}$, $\bar r=r/(1-r)$, and $\boldsymbol{x^\star}{=}\arg\min_{\boldsymbol{x}} f(\boldsymbol x)$. 
\end{prop}
\balance
\nocite{beck2017first}
\AtNextBibliography{\small} 
\printbibliography
% {\color{blue!80!black} \printbibliography[category=important, heading=none]}

\end{document}